\documentclass[11pt]{article}
\usepackage[margin=1in]{geometry}
\usepackage{amsthm,amsmath}
\usepackage{latexsym,amssymb}
\usepackage{graphicx,float,color,fancybox,shapepar,setspace,hyperref}
\usepackage{subfigure}
\usepackage{pgf,tikz}
\usepackage[affil-it]{authblk}
\usepackage{indentfirst}
\usepackage[section]{placeins}
\hypersetup
{
colorlinks=true,
linkcolor=blue,
filecolor=blue,
urlcolor=blue,
citecolor=cyan,
}
\usetikzlibrary{arrows}
\makeatletter
\renewenvironment{proof}[1][\proofname]{\par
  \normalfont \topsep6\p@\@plus6\p@\relax
  \trivlist
  \item[\hskip\labelsep
        \itshape
    #1\@addpunct{.}]\ignorespaces
}{%
  \hfill $\square$ 
  \endtrivlist\@endpefalse
}
\makeatother

\newtheorem{thm}{Theorem}
\newtheorem{conj}{Conjecture}

\newtheorem{pr}{Problem}
\newtheorem{claim}{Claim}[section]
\newtheorem{cor}{Corollary}
\newtheorem{lem}{Lemma}

 \def\qed{\hfill\square}

\def\~{\sim}

 \def\qed{ \hfill $\blacksquare$}

\begin{document}

\title{Fan-goodness of sparse graphs}
\author{Ting HUANG$^1$, Yanbo ZHANG$^2$, Yaojun CHEN$^{1,}$\footnote{Corresponding author. Email: yaojunc@nju.edu.cn}\\
{\small $^1$School of Mathematics, Nanjing University, Nanjing 210093, China}\\
{\small $^2$School of Mathematical Sciences, Hebei Normal University, Shijiazhuang 050024, China}}
 \date{}
\maketitle

\begin{abstract}
Let $G$ be a connected graph of order $n$, $F_k$ be a fan consisting of $k$ triangles sharing a common vertex, and $tF_k$ be $t$ vertex-disjoint copies of $F_k$. Brennan (2017) showed the Ramsey number $r(G,F_k)=2n-1$ for $G$ being a unicyclic graph for $n \geq k^2-k+1$ and $k\ge 18$, and asked the threshold $c(n)$ for which $r(G,F_k) \geq 2n$ holds for any $G$ containing at least $c(n)$ cycles and $n$ being large. In this paper, we consider fan-goodness of general sparse graphs and show that if $G$ has at most $n(1+\epsilon(k))$ edges, where $\epsilon(k)$ is a constant depending on $k$, then 
$$r(G,F_k)=2n-1$$ 
for $n\ge 36k^4$, which implies that $c(n)$ is greater than $\epsilon(k) n$. Moreover, if $G$ has at most $n(1+\epsilon(k,t))$ edges, where $\epsilon(k,t)$ is a constant depending on $k,t$, then
$$r(G,tF_k)=2n+t-2$$ 
provided $n\ge 161t^2k^4$.


\vskip 2mm
\noindent{\bf Keywords}: Ramsey number, sparse graph, fan
\end{abstract}

\section{Introduction}
Given two graphs $G$ and $H$, the \emph{Ramsey number} $r(G, H)$ is defined as the smallest positive integer $N$ such that, for every red-blue edge coloring of the complete graph $K_N$, there exists either a red subgraph isomorphic to $G$ or a blue subgraph isomorphic to $H$. When both $G$ and $H$ are complete graphs, $r(G, H)$ is the \emph{classical Ramsey number}. Let $H$ be a graph with chromatic number $\chi (H)$.
The \emph{chromatic surplus} $s(H)$ is defined to be the smallest number of vertices in a color class under any  $\chi (H)$-coloring of $V(H)$. For a connected graph $G$ with $n \geq s(H)$ vertices, Burr~\cite{Burr1981} established the following general lower bound for $r(G,H)$: 
\begin{equation} \label{lower}
        r(G, H)\ge (\chi(H)-1)(n-1)+s(H)\,.
    \end{equation}

Burr~\cite{Burr1981} also defined $G$ to be  \emph{$H$-good} if equality holds in inequality~(\ref{lower}). 
Prior to the formal definition of Ramsey goodness, Chv\'atal \cite{Chvatal1977} had established that all $n$-vertex trees $T_n$ satisfy the \emph{$K_k$-goodness} property, as demonstrated in the following result.
 \begin{thm}[Chv\'atal~\cite{Chvatal1977}]\label{Chv}  $r(T_n,K_k)=(n-1)(k-1)+1\,.$           \end{thm}

 For two graphs $G$ and $H$, 
 $G\cup H$ denotes their disjoint union, and 
 $G + H$ represents the disjoint union of two
 graphs $G$ and $H$, together with the new edges connecting every vertex of $G$  to every vertex of $H$. Moreover, $kH$ denotes the disjoint union of $k$ copies of $H$. A fan $F_k$ is a graph consisting of $k$ triangles with one vertex in common, i.e., $F_k=K_1+kK_2$.

 The complete-fan Ramsey number has attracted much attention.
 Because $F_1=K_3$, the earliest results on complete-fan Ramsey numbers follow from classical Ramsey theory. It is known \cite{AKS1980,AKS1981,S1983} that $r(K_n,K_3)\leq(1+o(1))n^2/\log n$. Kim \cite{Kim1995} established a lower bound of the correct order in 1995. The best currently available lower bound is due to Campos, Jenssen, Michelen, and Sahasrabudhe \cite{CJMS2025}:
 $$r(K_n,K_3)\geq\left(\frac13+o(1)\right)\frac{n^2}{\log n}.$$
 It follows that
$$r(K_n,K_3) = r(K_n,F_1)=\Theta(\frac{n^2}{\log n}).$$
 In 1996, Li and Rousseau \cite{LR1996} established an upper bound for $r(K_n,F_k)$, showing 
 $$r(K_n,F_k) \leq (1+o(1))\frac{n^2}{\log n}.$$
Moreover, Li and Rousseau \cite{LR1996}  conjectured that $r(K_n,F_1)=(1+o(1))n^2 / \log n$, and by the fact that $r(K_n,F_1)\leq r(K_n,F_k)\leq (1+o(1))n^2 / \log n$, they proposed a more general conjecture as below. 
\begin{conj} (Li and Rousseau \cite{LR1996})
If $k$ is any fixed positive integer, then
    $$r(K_n,F_k) = (1+o(1))\frac{n^2}{\log n} ~(n \to \infty).$$
\end{conj}
This conjecture is still open. However, for cases where the complete graph is sufficiently small relative to the fan, several results are available, see \cite{CL2025,KOS2019,LR1996,BB2005,WQ2022,ZC2023}.



Subsequent studies have extended this to cases where 
$K_n$ is replaced by non-complete graphs (e.g., cycles or trees), which have garnered increasing interest.
In 2010, Shi \cite{Shi 2010} proved that $C_n$ is $F_k$-good for $n \geq 3k+1$.
 In 2015, Zhang, Broersma, and Chen \cite{ZBC} showed that the star $K_{1,n-1}$ is $F_k$-good and $T_n$ is $F_k$-good as below.
\begin{thm} (Zhang, Broersma, and Chen \cite{ZBC}) \label{starfan}
 $r(K_{1,n-1},F_k) = 2n-1$ for $n\geq k^2 -k + 1$ and $k \neq 3,4,5$, and the lower bound $n\geq k^2 -k + 1$ is best
 possible. $r(K_{1,n-1},F_k) = 2n-1$ for $n \geq 6(k- 1)$ and $k= 3,4,5$.
\end{thm}
\begin{thm}(Zhang, Broersma, and Chen \cite{ZBC}) \label{treefan}
$r(T_n, F_k) = 2n-1$ for all integers $k$ and $n \geq 3k^2 -2k -1$.
\end{thm}
Moreover, Zhang, Broersma and Chen \cite{ZBC}
 conjectured that $n \geq k^2 -k +1$ is the best achievable lower bound
 on $n$ in terms of $k$ for which $T_n$ is $F_k$-good. Brennan \cite{Brennan} confirmed this conjecture  for the case $k\geq 9$.
Brennan \cite{Brennan} also extended $F_k$-goodness from $T_n$ to unicyclic graphs $UC_n$, where $UC_n$ denotes 
a connected graph with $n$ vertices and $n$ edges.
\begin{thm}(Brennan \cite{Brennan}) \label{Brennan}
If $k \geq 18$ and $n \geq k^2 -k + 1$, then
$r(UC_n, F_k) = 2n-1$.
\end{thm}
Note that trees contain no cycle and unicyclic graphs contain exactly one cycle and they are both $F_k$-good.
 Brennan proposed an interesting problem: 
 \begin{pr} (Brennan \cite{Brennan})
 Given an integer $k$, what is the threshold $c(n)$ such that if $G$ is any connected graph with $n$ vertices and at least $c(n)$ cycles, then $r(G,F_k) \geq 2n$ (i.e., $G$ is not $F_k$-good)?   
 \end{pr}

 Brennan's result implied that if $k \geq 18$ and $n \geq k^2-k+1$, then $c(n)\geq 2$. 
Obviously, any edge addition to a connected graph produces at least a new cycle, and the number of cycles increasing monotonically with the number of added edges.
 Brennan's problem motivates the study of fan-goodness of sparse graphs. The following theorem implies that, for each fixed positive integer $k$, $c(n)$ grows at least linearly with $n$.


\begin{thm} \label{mr}
 Let $k$ be a positive integer and $n\geq 36k^4 $. If $G$ is a connected graph with $n$ vertices
and at most $ n(1+1/(204k^3+126k^2))$ edges,
then
$$r(G,F_k)=2n-1.$$
\end{thm}

Indeed, adding $\lfloor n/(204k^3+126k^2)\rfloor$ distinct chords to $C_n$ gives a graph covered by Theorem~\ref{mr}. Relative to a spanning path of $C_n$, the added chords and the remaining cycle edge determine $\lfloor n/(204k^3+126k^2)\rfloor+1$ distinct fundamental cycles. Thus $c(n)>n/(204k^3+126k^2)$ for $n\geq36k^4$.

Furthermore, if $F_k$ is replaced by the disjoint union of copies of $F_k$, an intriguing question arises: which sparse graphs are still $tF_k$-good? In this paper, we also investigate this problem and prove the following two results.
\begin{thm}\label{star tfan}
  Let $k,t$ be positive integers and $n\geq \max\{12tk+2k,6tk^2\}$. Then $$r(K_{1,n-1},tF_k)=2n+t-2.$$
\end{thm}

\begin{thm} \label{mr2}
  Let $k,t$ be positive integers and $n\geq 161t^2k^4$. If $G$ is a connected graph with $n$ vertices
and at most $ n(1+1/(204 t k^3 + 147 t k^2 ))$ edges, 
then
$$r(G,tF_k)=2n+t-2.$$
\end{thm}
 
The organization of this paper proceeds as follows. In Section~\ref{section2}, we introduce some basic lemmas needed for the proofs. In Sections~\ref{section3}, \ref{section4}, and \ref{section5}, we present the proofs of Theorems \ref{mr}, \ref{star tfan}, and \ref{mr2}, respectively.

We conclude this section by introducing some additional notation. Throughout this paper, all graphs are finite and simple without loops.   
 For a graph $G=(V(G),E(G))$, let $e(G)=|E(G)|$. 
For $v\in V(G)$, the neighborhood of $v$ is $N_G(v)=\{u~|~uv\in E(G)\}$ and $d(v)=|N_G(v)|$. The maximum degree of $G$ is $\Delta(G)=\max\{d(v)~|~v\in V(G)\}$.
The
 graph $G-H$ will be the subgraph of $G$ induced by the vertices of $G$ not in $H$. Given a red-blue edge-coloring of $K_N$, $K_N[R]$ and $K_N[B]$ denote the spanning subgraphs formed by the red edges and the blue edges, respectively.
 For a vertex $v$, its red neighborhood $N_R(v)$ (resp.\ blue neighborhood $N_B(v)$) consists of vertices adjacent via red (resp.\ blue) edges, with corresponding degrees $d_R(v) = |N_R(v)|$ and $d_B(v) = |N_B(v)|$. The maximum blue degree of $G$ is $\Delta_B(G)=\max\{d_B(v)~|~v\in V(G)\}$.
 For a complete graph $K_N$ and a vertex subset $A$, we use $K_N[A]$ to denote
 the subgraph induced by $A$. The independence number of $G$ is denoted by $\alpha(G)$. 
 A path $P$ of $G$ is a suspended path if each vertex of $P$, except for its endvertices, has degree $2$ in $G$. An end-edge is one incident with a vertex of degree 1.

\section{Some Basic Lemmas}\label{section2}

In 1982, Burr, Erd\H{o}s, Faudree, Rousseau, and Schelp~\cite{BEFRS} showed that any sparse graph must contain a long suspended path 
or many degree-1 vertices, thereby characterizing a useful structural feature of sparse graphs. We use the following form of their result, with cycle components excluded.
\begin{lem} (Burr, Erd\H{o}s,  Faudree, Rousseau, and Schelp \cite{BEFRS})  \label{old}
 Let $G$ be a graph on $n$ vertices and $n+\ell$ edges. If $G$ has no isolated
vertices, no component that is a cycle, and no suspended path with more than $q$ vertices, then $G$ has at least
$\lceil \frac{n}{2q}-\frac{3\ell}{2}\rceil$ vertices of degree 1.   
\end{lem} 
Zhang and Chen \cite{ZC} recently developed an enhanced version of the Trichotomy Lemma for sparse graphs, which provides the clearer structure of a sparse graph. 
\begin{lem}(Zhang and Chen \cite{ZC}) \label{trichotomy}
Let $G$ be a connected graph with $n$ vertices and $n+\ell$ edges that is not a star, where $\ell \geq -1$, $n\geq q \geq 3$, and $s\geq2$ is an integer. If $G$ contains neither a suspended path of order $q$ nor a matching consisting of $s$ end-edges, then the number of vertices of degree at least 2 in $G$ is at most $\gamma$, and 
$G$ has a vertex adjacent to at least $\lceil \frac{n-\gamma}{s-1} \rceil$ vertices of degree 1, where $\gamma=(q-2)(2s+3\ell-2)+1$.
\end{lem}
 By using the Trichotomy Lemma for sparse graphs as above, we divide the problem into three cases. In each case, we use Lemma \ref{weak} to identify a subgraph $H$ of $G$, where $G$ is the given sparse graph. Two tools are utilized to extend $H$ into the desired graph $G$. The first case relies on the following path extension lemma in the form given by Erd\H{o}s, Faudree, Rousseau, and Schelp~\cite{EFRS1985}, based on an argument of Bondy and Erd\H{o}s~\cite{BE}.
 For the second case, Hall's theorem is applied to ensure the embedding of a matching. 
 In the third case, it is necessary to first identify a star $K_{1,n-1}$, where Theorems \ref{starfan} and \ref{star tfan} play a pivotal role.

 \begin{lem} (Erd\H{o}s, Faudree, Rousseau, and Schelp \cite{EFRS1985})\label{path extension}
    Let $K_{a+b}$ be a complete graph on the vertex set
 $\{x_1,\ldots,x_a,y_1,\ldots,y_b\}$, with edges colored red or blue. Assume that there is a red path
 $x_1x_2\cdots x_a$ of length $a-1$ joining $x_1$ and $x_a$. If $a 
\geq b(c-1)+d$, then at least one of
 the following holds:
 
 \noindent
 (1) A red path of length a connects $x_1$ to $x_a$;
 
 \noindent
 (2) A blue complete subgraph $K_c$ exists;
 
 \noindent
 (3) There are $d$ vertices in $\{x_1,x_2,\ldots,x_a\}$ that are joined in blue to every vertex in
 $\{y_1,y_2,\ldots,y_b\}$.
     
 \end{lem}

\begin{lem} (Hall \cite{Hall}) \label{matching}
 Consider a complete bipartite graph $K_{a,b}$, where $a \leq b$, with parts $X = \{x_1,\ldots,x_a\}$ and $Y=\{y_1,y_2,\ldots,y_b\}$, whose edges are colored red and blue. Then, one of the following holds:

 \noindent
 (1) There exists a red matching of size $a$;

 \noindent
 (2) For some $0 \leq c \leq a-1$, there exists a blue subgraph $K_{c+1,b-c}$, where $c+1$ vertices are in $X$.   
 \end{lem}

 Moreover, the proofs of Theorems~\ref{mr} and~\ref{mr2} require an upper bound for the Ramsey number of any graph $G$ versus $F_k$ in terms of the order and size of $G$, which is interesting in its own right. First, we need the following upper bound for $r(G,kK_2)$.
\begin{lem} (Faudree, Schelp, and Sheehan \cite{FSS}) \label{FSS}
Let $G$ be a simple graph on $n$ vertices none of which is isolated and $\alpha(G)=\alpha$. For any positive integer $k$,
$$r(G,kK_2) \leq \max\{n+2k-2-\left\lfloor \frac{2\alpha}{3} \right\rfloor, n+k-1\}.$$
    
\end{lem}

\begin{lem} \label{upper}
    If $G$ is a simple graph on $n\geq1$ vertices with $m$ edges, then for every positive integer $k$,
     $$r(G,F_k) \leq n+2mk-\frac{2m}{n} .$$
\end{lem}
\noindent
{\bf Proof.}
We use induction on $n$. The assertion is immediate for $n=1$. Let $v$ be a vertex of minimum degree $\delta$. If $\delta=0$, then
\begin{align*}
r(G,F_k)\le r(G-v,F_k)+1 \le n+2mk-\frac{2m}{n-1}
\le n+2mk-\frac{2m}{n},
\end{align*}
so assume $\delta\ge1$. Put
$N=\max\{r(G-v,F_k),\,\delta(r(G,kK_2)-1)+n\}$.
Assume that $K_N$ has no blue $F_k$. It contains a red $G-v$. Let $A$ be the images of the neighbors of $v$ in this copy, and let $B$ be the vertices outside the copy. If some vertex of $B$ is red-adjacent to all of $A$, it extends the copy to a red $G$. Otherwise every vertex of $B$ sends a blue edge to $A$. Since
$|B|=N-n+1\ge\delta(r(G,kK_2)-1)+1$,
some $w\in A$ has at least $r(G,kK_2)$ blue neighbors. Their induced coloring contains either a red $G$ or a blue $kK_2$; the latter, together with $w$, is a blue $F_k$. Hence $r(G,F_k)\le N$.

As $G$ has no isolated vertices, Lemma~\ref{FSS} gives $r(G,kK_2)\le n+2k-2$. The induction hypothesis therefore yields
$$r(G,F_k)\le\max\left\{n-1+2k(m-\delta)-\frac{2(m-\delta)}{n-1},\ n+\delta(n+2k-3)\right\}.$$
For the first term, using $2m\ge n\delta$, we have
\begin{align*}
&n+2km-\frac{2m}{n}-\left(n-1+2k(m-\delta)-\frac{2(m-\delta)}{n-1}\right)\\
&\qquad=1+2k\delta+\frac{2(m-n\delta)}{n(n-1)}
\ge1+2k\delta-\frac{\delta}{n-1}>0.
\end{align*}
For the second term, $\delta\le2m/n$ gives
\begin{align*}
n+\delta(n+2k-3)\le n+2m+\frac{4km}{n}-\frac{6m}{n} \le n+2km-\frac{2m}{n},
\end{align*}
where the last difference is $2m(k-1)(1-2/n)\ge0$. This proves the assertion. \qed

\vskip 2mm
Because $F_k=K_1+kK_2$, identifying $F_k$ in a graph can be reduced to finding a $kK_2$ within the neighborhood of some vertex. This necessitates the following two lemmas. 
\begin{lem} (Chv\'atal and Harary \cite{Chvatal1972}) \label{K2}
For any graph $G$ on $n$ vertices that contains no
 isolated vertices and $G$ is not complete, 
 $r(G,2K_2)=n+1$.
    
\end{lem}

\begin{lem} (Zhang and Chen \cite{ZC}) \label{kK2}
Let $k\geq 3$ be an integer and $n \geq 3(4k-5)(2k-3)$.
     If $G$ is a connected graph with $n$ vertices and  at most $n+n^2/(4k-5)-2$ edges, then
 $r(G,kK_2) = n+k-1$.
    
\end{lem}
\noindent Moreover, it is trivial that $r(G,K_2)=n$ for any graph $G$ on $n$ vertices.

\section{Proof of Theorem \ref{mr}} \label{section3}
Note that $F_k\subseteq \overline{K}_k+K_{1,k}$, the following lemma can help us first identify a blue $K_{1,k}$ and then find a blue $F_k$. 
\begin{lem} (Huang, Zhang, and Chen \cite{HZC}) \label{n+k-1}
For integers $k \geq 1$ and $n\geq 6k^3$, let $G$ be a connected graph with $n$ vertices and at most $ n(1+1/(24k-12))$ edges. Then
$$r(G,K_{1,k})\leq n+k-1.$$
\end{lem}

\begin{lem} \label{weak}
    For integers $k \geq 1$ and $ n\geq 6k^3$, let $G$ be a connected graph with $n$ vertices
and at most $ n(1+1/(6k^3+12k^2+11k/2))$ edges. Then $$r(G,F_k)\leq 2n+k-2.$$ 
\end{lem}

\noindent
{\bf Proof.}
Let $N=2n+k-2$. Assume that $K_N$ is a complete graph with a red-blue
 edge-coloring such that it contains neither a red $G$ nor a blue
  $F_k$. 

 First, let $G_1\subseteq G$ be obtained from $G$ by recursively deleting vertices of degree 1, one at a time, until the resulting graph has no vertex of degree 1. We prove, by reversing these deletions, that $K_N$ contains no red $G_1$. Let $G_1 \subseteq G'\subseteq G$, and let $G''$ be obtained from $G'$ by deleting a vertex of degree 1. Assume that $K_N$ has no red $G'$ and, to the contrary, contains a red $G''$. If $v$ is the image of the support vertex of the deleted leaf, every vertex outside this red copy of $G''$ is blue-adjacent to $v$; otherwise the copy extends to a red $G'$. 
For $k\geq3$, $|G| \geq 6k^3\geq 3(4k-5)(2k-3)$ and 
$$e(G)\leq n\left(1+\frac{1}{6k^3+12k^2+11k/2}\right)\leq n+\frac{n^2}{4k-5}-2,$$
 so Lemma~\ref{kK2} gives $r(G,kK_2)=n+k-1$. For $k=2$, Lemma~\ref{K2} gives the same equality, and for $k=1$ the equality $r(G,K_2)=n$ is trivial. Consequently,
 $$|N_B(v)|\geq|K_N-G''|\geq N-(n-1)=n+k-1 \geq r(G,kK_2).$$  
 Thus $K_N[N_B(v)]$ contains either a red $G$ or a blue $kK_2$. The first alternative contradicts the global assumption, while the second, together with $v$,
 contains a blue $K_1+kK_2=F_k$, a contradiction.

 Secondly, let $G_0$ be the graph obtained from $G_1$ by shortening any suspended path with at least $(2k+3)k+1$ vertices to one with $(2k+3)k$ vertices.
 In this process, let $G''$ be a graph obtained from $G'$ by shortening any suspended path with at least $(2k+3)k+1$ vertices by a vertex.
We will prove that $K_N$ contains no red $G_0$
 by showing that if $K_N$ has no red $G'$, then $K_N$ contains no red $G''$ either. Suppose to the contrary that $K_N$ contains a red $G''$. Because $|G| \geq 6k^3$ and 
\begin{align*}
     e(G) &\leq  n\left(1+\frac{1}{6k^3+12k^2+11k/2}\right) \leq n\left(1+\frac{1}{24k-12}\right),
     \end{align*}  
     by Lemma \ref{n+k-1}, we have $|K_N-G''|\geq N-(n-1)=n+k-1 \geq r(G,K_{1,k})$. This implies that $K_N-G''$ contains a blue $K_{1,k}$. Let $P=x_1x_2\cdots x_a$ be the suspended path in $G''$, where $a\geq (2k+3)k=(k+1)[(2k+1)-1]+k$. Let the vertices of the blue star be relabeled as $y_1,y_2,\ldots,y_b$, where $b=k+1$. Also, let $c=2k+1$ and $d=k$. So we have $a \geq b(c -1)+d$. The first alternative of Lemma~\ref{path extension} would extend the red $G''$ to a red $G'$, contrary to the choice of $G'$. Hence, by the lemma, we either have a blue subgraph $K_{2k+1}$ or
there exist $k$ vertices among $\{x_1,x_2,\ldots,x_a\}$ such that every edge between these $k$ vertices and the vertices $\{y_1,y_2,\ldots,y_b\}$ is blue. In the first case, we obtain a blue $K_{2k+1}$,
 and in the second case, we obtain a blue $\overline{K}_k+K_{1,k}$. In both cases, there
 exists a blue $F_k$, again a contradiction.

Assume $G_0$ has $\ell$ vertices and $m$ edges, then $m\leq \ell+n/(6k^3+12k^2+11k/2)$. Moreover, since $K_N$ contains no red $G_0$ and $G_0$ is connected, we must have $\ell \geq 3$ and $m\geq \ell-1$. Recall that $G_0$ has no vertices of degree 1 and no suspended path with more than $(2k+3)k$ vertices.
If $G_0=C_\ell$, then $\ell\leq(2k+3)k$. For $k\geq2$, Lemma~\ref{upper} gives
\[
r(C_\ell,F_k)\leq(2k+1)\ell-2
\leq(2k+1)(2k+3)k-2
\leq12k^3+k-2\leq N.
\]
For $k=1$, we have $\ell\in\{3,4,5\}$ and $N=2n-1\geq11$; the classical equality $r(C_3,K_3)=6$ and Shi's result~\cite{Shi 2010}, which gives $r(C_\ell,K_3)=2\ell-1$ for $\ell=4,5$, again yield $r(C_\ell,F_1)\leq N$. Both cases contradict the absence of a red $G_0$. Thus $G_0$ is not a cycle, and Lemma \ref{old} applies to give
$$\frac{\ell}{2(2k+3)k}-\frac{3n}{2(6k^3+12k^2+11k/2)}\leq 0.$$
Therefore, $\ell \leq \frac{3(2k+3)k}{6k^3+12k^2+11k/2}n$.  Hence, by Lemma \ref{upper},
\begin{align*}
    r(G_0,F_k)\leq &\ell+2km-\frac{2m}{\ell}\\
    \leq& \ell +2k\left(\ell+\frac{n}{6k^3+12k^2+11k/2}\right)-\frac{2(\ell-1)}{\ell}\\
    \leq & \ell(2k+1)+\frac{2k}{6k^3+12k^2+11k/2}n-1 \\
    \leq
    & \frac{12k^3+24k^2+11k}{6k^3+12k^2+11k/2}n-1\\
    =&2n-1
    \leq N.
\end{align*}
This is to say that $K_N$ has a red $G_0$, a contradiction, and this completes the proof. \qed

\vskip 5mm
We now begin to  prove Theorem \ref{mr}.

Since $\chi(F_k)=3$ and $s(F_k)=1$, \eqref{lower} gives the matching lower bound. It therefore suffices to prove the upper bound. Let $N=2n-1$. We proceed by contradiction. Assume that $K_N$ is a complete graph with a red-blue
 edge-coloring such that it contains neither a red $G$ nor a blue $F_k$. We will derive a contradiction in the following argument.

 If $k=1$, the assertion follows directly from Lemma~\ref{weak}, since the hypotheses of Theorem~\ref{mr} imply those of that lemma. Hence assume $k\ge2$. By Theorem \ref{starfan}, the case where $G$ is a star is already settled. Thus we assume $G$ is not a star, and then 
 by Lemma \ref{trichotomy}, we distinguish three cases.

\vskip 2mm 
\noindent 
{\textbf{Case 1.}} $G$ has a suspended path $P$ with $(2k+3)k+k-1$ vertices.
\vskip 2mm

 Let $H$ be the graph obtained from $G$ by shortening $P$ by $k-1$ vertices. Since $n\geq 36k^4$, we have 
 \begin{align*}
 |H| &=n-k+1 \geq  6k^3 ,\\
   e(H) &\leq n\left(1+\frac{1}{204k^3+126k^2}\right)-(k-1)  
\leq \left(n-k+1 \right) \left(1+\frac{1}{6k^3+12k^2+11k/2}\right).
 \end{align*}
 Then by Lemma \ref{weak}, $$r(H,F_k)\leq 2(n-k+1)+k-2=2n-k\leq 2n-1=N.$$ 
 Hence $K_N$ contains a red $H$. Let $H'$ be the subgraph of $K_N[R]$ which can
 be obtained from $H$ by lengthening the suspended path as much as possible (up to
 $k-1$ vertices). If $H'$ is $G$, the proof of this case is complete. If not, there are two possibilities to consider, according to whether $K_N-H'$ has a blue $K_{1,k}$ or not.
 
 If $K_N-H'$ contains a blue $K_{1,k}$, let $P'= x_1 x_2 \cdots x_a$ be the suspended path in $H'$, where $a\geq (2k+3)k=(k+1)((2k+1)-1)+k$. Let $y_1,y_2,\ldots,y_b$ be the vertices of the blue $K_{1,k}$, where $b=k+1$. Set $c=2k+1$ and $d=k$.  So we have $a \geq b(c -1)+d$. Applying Lemma \ref{path extension} on the red path $P'$ and $\{y_1,y_2,...,y_b\}$, its first alternative contradicts the maximality of $P'$. Hence either there is a blue $K_{2k+1}$ or there are $k$ vertices among $\{x_1,x_2,\ldots,x_a\}$ such that every edge between these vertices and $\{y_1,y_2,\ldots,y_b\}$ is blue. The first outcome contains a blue $F_k$, and the second contains a blue $\overline{K}_k+K_{1,k}$ and hence a blue $F_k$, a contradiction.

 If $K_N-H'$ has no blue $K_{1,k}$, then $k\geq 2$. Let  $S=V(K_N-H')$. 
Because $ |H'| \le n-1$, we have $ |S|\geq N-(n-1)= n$. 
Since $n\ge 36k^4$, we have \begin{align*}  |H| &=n-k+1 \geq  6k^3,\\   
e(H) &\leq n\left(1+\frac{1}{204k^3+126k^2}\right)-(k-1)  
\leq \left(n-k+1 \right) \left(1+\frac{1}{24k-12}\right).  \end{align*}         
Then by Lemma \ref{n+k-1}, $$r(H,K_{1,k})\leq n-(k-1)+k-1=n\le |S|.$$ Because $K_N[S]$ contains no blue $K_{1,k}$, there must exist a red $H$ in $K_N[S]$. 
Let $H''$ be the subgraph of $K_N[S]$ obtained from $H$ by lengthening the shortened suspended path through red edges as many as possible (up to
$k-1$ vertices). Since $K_N$ has no red $G$, $H''\not=G$. 
Let $P''$ be the suspended path of $H''$. Since 
$|S|-|H''| \geq n-(n-1)=1$, we may assume $v\in S\backslash V(H'')$. 
By the maximality of $P''$, $v$ is not red-adjacent to two consecutive vertices of $V(P'')$ in $K_N$, and so $v$ has at least $\lfloor |V(P'')|/2 \rfloor  \geq \lfloor(2k+3)k/2\rfloor \ge k$ blue neighbors in $V(P'')$. Thus, there must exist a blue $K_{1,k}$ with center $v$ in $K_N[S]$, a contradiction.


 \vskip 2mm
  \noindent  {\textbf{Case 2.}} $G$ has a matching consisting of $2k-2$ end-edges.
 \vskip 2mm
Let $H$ be the graph obtained from $G$ by deleting $2k-2$ end-vertices of this matching. Since $n\geq 36k^4$, we have 
\begin{align*}
 |H| &=n-(2k-2) \geq 6k^3,\\
   e(H) &\leq n \left(1+\frac{1}{204k^3+126k^2}\right)-(2k-2)  \\
   &\leq (n-(2k-2))\left(1+\frac{1}{6k^3+12k^2+11k/2}\right).
 \end{align*}
By Lemma \ref{weak}, 
$$r(H,F_k)\leq 2(n-(2k-2))+k-2 =2n-3k+2\leq 2n-1= N.$$ 
Thus $K_N$ contains a red $H$. 
Let $X$ be the set of $2k-2$ vertices of $H$ which were adjacent to the end-vertices deleted from $G$ and $Y=V(K_N)-V(H)$. By Lemma \ref{matching}, there is either a red matching that covers $X$, or a blue $K_{c+1,|Y|-c}$, where $0 \leq c \leq |X|-1$. 
If the former holds, then $K_N$ would contain a red $G$, a contradiction. If the latter holds, let $X'$ be the vertices of $X$ and $Y'$ be the vertices of $Y$ in the blue $K_{c+1,|Y|-c}$. 
 Since $K_N$ contains no blue $F_k$, there is no blue $kK_2$ in  $K_N[N_B(v)]$ for any $v\in V(K_N)$, and then
 $|N_B(v)|\leq n+k-2$ by Lemmas \ref{K2} and  \ref{kK2}. Thus  $|Y'|=|Y|-c\leq n+k-2$.
 Then \begin{align*}
     c&\geq |Y|-n-k+2\\
     &=2n-1-(n-|X|)-n-k+2\\
     &= |X|-k+1\\
     &= k-1.
     \end{align*}
     Moreover, 
 $|Y'|=|Y|-c\geq 2n-1-(n-|X|)-(|X|-1)=n$.
 
 If $K_N[Y']$ has a blue $kK_2$, then there is a blue $F_k$ with center in $X'$, and if $K_N[Y']$ has a blue $K_{1,k}$, then since $|X'|=c+1\geq k$, there is a blue $\overline{K}_k+K_{1,k}$, which contains an $F_k$. Thus $K_N[Y']$ contains neither a blue $kK_2$ nor a blue $K_{1,k}$. 
 
Select $Y''\subseteq Y'$ with $|Y''|=n$. Let $M$ be a maximum blue matching in $K_N[Y'']$, and write $s=|V(M)|$. Since there is no blue $kK_2$, $s\leq2k-2$, and $K_N[Y''-V(M)]$ is red complete. Moreover, $|N_B(v)\cap (Y''-V(M))|\leq k-2$ for each $v\in V(M)$, because $v$ is blue-adjacent to its partner in $M$ and $K_N[Y']$ contains no blue $K_{1,k}$. In particular,
$$|N_R(v)\cap (Y''-V(M))|\geq n-s-(k-2)>s,$$
so a greedy choice gives a red matching from $V(M)$ into $Y''-V(M)$ covering $V(M)$. Delete $s$ end-vertices from the fixed matching of $2k-2$ end-edges of $G$, obtaining a graph $H'$ of order $n-s$. First choose the covering red matching, and then embed $H'$ into the red clique $Y''-V(M)$ so that the $s$ distinct support vertices are mapped to the corresponding matched vertices. Mapping the deleted leaves to $V(M)$ now gives a red $G$, a contradiction.


 \vskip 2mm
\noindent {\textbf{Case 3.}} 
$G$ has neither a suspended path of $(2k+3)k+k-1$ vertices, nor a matching formed by $2k-2$ end-edges.
\vskip 2mm

To apply Lemma \ref{trichotomy} in a more straightforward manner, we let $e(G)=n+\ell$, $q=(2k+3)k+k-1=2k^2+4k-1$, and $s=2k-2$, where $\ell\leq n/(204k^3+126k^2)$.

Let $G'$ be a graph obtained from $G$ by deleting all vertices of degree 1. Clearly, $G'$ is connected. Because $n\geq 36k^4$,
 then by Lemma \ref{trichotomy}, we have
\begin{align*}
|G'| &\leq \gamma=(q-2)(2s+3\ell-2)+1\\  
&\leq \left( 2k^2+4k-3\right) \left(4k-6+\frac{3n}{204k^3+126k^2}\right)+1 \\   &= 8k^3+4k^2-36k+19+\frac{3(2k^2+4k-3)}{204k^3+126k^2}n\\ 
&\leq \frac{96k^2}{204k^3+126k^2}n+\frac{3(2k^2+4k-3)}{204k^3+126k^2}n\\ 
& = \frac{102k^2+12k-9}{204k^3+126k^2}n.
\end{align*}
Moreover,
\begin{align*}
|G'|-1 \leq e(G') &\leq |G'|+\frac{n}{204k^3+126k^2}\leq \frac{102k^2+12k-8}{204k^3+126k^2}n.
 \end{align*}
By Lemma \ref{trichotomy},  $G$ has a vertex adjacent to at least 
$$\left\lceil \frac{n-\gamma}{2k-3} \right\rceil$$ 
vertices of degree 1.
Denote this vertex by $v$. By Theorem \ref{starfan}, there exists a red star $K_{1,n-1}$ in $K_N$, with the center $x$. By Lemma \ref{upper}, if $|G'| \geq 2$, then
\begin{align*}
    r(G',F_k)&\leq |G'|+2ke(G')-\frac{2e(G')}{|G'|}\\
    &\leq  \frac{102k^2+12k-9}{204k^3+126k^2}n+\frac{2k(102k^2+12k-8)}{204k^3+126k^2}n -\frac{2(|G'|-1)}{|G'|}\\
    &\leq \frac{204k^3+126k^2-4k-9}{204k^3+126k^2}n-1 \\
    &< n-1.
\end{align*}
Thus $K_N[N_R(x)]$ contains a red $G'$. Moreover, this holds trivially if $|G'| = 1$. Replacing $v$ in this $G'$ by $x$ gives another copy of red $G'$. Every vertex of $G''-G'$ is an original leaf of $G$ whose support lies in $G'$. 
Next, based
on the embedding of $G'$,
we apply a greedy algorithm to
embed $G''$ into
 the red subgraph of $K_N$, where $G''$ is obtained by removing all degree-1 vertices adjacent to $v$ in the graph $G$. Clearly, $G''$ contains $G'$ as a subgraph and
 $$|G''| \leq n-\left\lceil \frac{n-\gamma}{2k-3}\right\rceil.$$
If at some step the embedding cannot proceed, the image $z$ of the relevant support vertex has no unused red neighbor and hence is incident with at most $|G''|-2$ red edges. In other words, the
 number of blue edges incident with $z$ is at least $N-|G''|+1$.
Because $n\geq 36k^4$,
we have $n-\gamma \geq (k-1)(2k-3)$.
Then by Lemmas \ref{K2} and \ref{kK2}, we have $$|N_B(z)| \geq N-|G''|+1\geq 2n-1-n+\left\lceil \frac{n-\gamma}{2k-3}\right\rceil +1\geq n+k-1 =r(G,kK_2).$$ Then in the blue neighborhood of 
 $z$, there is a blue $kK_2$, which along with $z$ gives a blue $F_k$. This leads to a contradiction.

 Thus, $G''$ can always be embedded into the red subgraph of $K_N$. Recall the vertex $v$ in $G''$ is embedded into $x$, and $x$ is incident with at least $n-1$ red edges, we can get a red $G$ in $K_N$, a contradiction. This completes the proof of Theorem \ref{mr}. \qed

\section{Proof of Theorem \ref{star tfan}}\label{section4}

\begin{lem} (Hu and Peng \cite{HP2023})\label{star kK_2}
  If $n\geq k+1$, then $r(K_{1,n-1}, kK_2)=n+k-1$.
\end{lem}

Now, we begin to prove Theorem \ref{star tfan}. Since $\chi(tF_k)=3$ and $s(tF_k)=t$, \eqref{lower} gives the matching lower bound, so it suffices to establish the upper bound. 

The stated hypothesis implies
$$n\geq \max\{12kt+2k-19,\,10kt+2k-13,\,6k^2t-4k^2-k-t+1\}.$$
We use induction on $t$. If $t=1$, the assertion follows
 by Theorem \ref{starfan}. Assuming the theorem holds for $t -1$, we now proceed to consider the case where $k\geq 1$ and $t \geq 2$.

Let $N=2n+t-2$. Given a red-blue
 edge-coloring of  $K_N$, assume that it contains neither a red $K_{1,n-1}$ nor a blue $tF_k$. 
By the induction hypothesis, we have
$$r(K_{1,n-1},(t-1)F_k)= 2n+t-3<N.$$
Thus $K_N$ contains a blue $(t-1)F_k$.
Let $A$ be the
 vertex set of this $(t-1)F_k$ and $S=V(K_N)-A$.   
 Clearly, 
 \begin{equation*}
    |S|=2n+t-2-(2k+1)(t-1)=2n-2kt+2k-1. 
 \end{equation*}

The main idea for proving Theorem \ref{star tfan} is as follows.
Firstly, we identify two large disjoint subsets $Z_1$ and $Z_2$ in $S$, each of which induces a nearly red complete subgraph. Second, we partition $S\setminus (Z_1\cup Z_2)$ into $W_1$ and $W_2$, and show $S_i=Z_i\cup W_i$  induces nearly red complete subgraphs for $i=1,2$.
Third, we further partition the vertex set $A$, and show that there exists a subset of $A$ together with some $S_i$ contains a red $K_{1,n-1}$ as a subgraph. This will lead to a contradiction.
\vskip 2mm
Since $K_N$ has no red $K_{1,n-1}$, we have 
\begin{equation*}
    |N_B(v)| \geq N-1-(n-2)=n+t-1\text{ for any $v\in V(K_N)$}.
\end{equation*}

\begin{claim} \label{1}
  For 
each $v\in S$,  $n-2kt+2k\le |N_B(v)\cap S| \leq n+k-2$. 
\end{claim}
\begin{proof} Let $v$ be any vertex in $S$. On one hand, we have 
$|N_B(v)\cap S| \geq n+t-1-|A|=n-2kt+2k$.
On the other hand, since $K_N[N_B(v)\cap S]$ contains no red $K_{1,n-1}$, and $K_N$ has no blue $tF_k$ implies $K_N[N_B(v)\cap S]$ contains no  blue $kK_2$, we have
$|N_B(v)\cap S| \leq r(K_{1,n-1},kK_2)-1$. Thus, by Lemma \ref{star kK_2}, we can deduce 
$|N_B(v)\cap S| \leq n+k-2$,
and so the assertion follows.
\end{proof}

\begin{claim} \label{Z1Z2}
    There are two disjoint sets $Z_1$ and $Z_2$ in $S$  such that for each $i=1,2$,
the blue subgraph in $K_N[Z_i]$ contains at most $3(k-1)$ nonisolated vertices and has maximum degree at most $2(k-1)$,
    and $n-2kt+k+1\le |Z_i| \le n+k-2$.
\end{claim}
\begin{proof}
For any $v\in S$, let $M_v$ be a blue maximum matching in $K_N[N_B(v)\cap S]$ and $X_v =N_B(v)\cap S- V(M_v)$. Obviously $K_N[X_v ]$ contains no blue edges, and $|M_v|\leq k-1$, otherwise $K_N[S]$ contains a blue $F_k$, a contradiction. Moreover, by the
 maximality of $M_v$, for $xy \in M_v$, if $|N_B(x)\cap X_v|\geq 2$, then $|N_B(y)\cap X_v|=0$; and if $|N_B(x)\cap X_v|=|N_B(y)\cap X_v|=1$, then $x$ and $y$ are adjacent to the
 same vertex in $X_v$. Let $Y_v\subseteq V(M_v)$ be the set of vertices that have at least two blue neighbors in $X_v$, and let $Z_v= N_B(v)\cap S-Y_v$. It is
 obvious that $|Y_v|\leq k-1$ and then by Claim \ref{1}, we have $$|Z_v|\geq |N_B(v)\cap S|-(k-1)\geq n-2kt+2k-(k-1)=n-2kt+k+1,$$ 
 and
 $|Z_v|\leq |N_B(v)\cap S|\leq n+k-2$.
 Moreover, $K_N[Z_v]$ contains at most $3(k-1)$ vertices whose blue degree is not equal to $0$ and the
 maximum blue degree of $K_N[Z_v]$ is at most $2(k-1)$.
 
 Since $X_v\subseteq Z_v$ and $|X_v| \geq n-2kt+2k- 2(k-1)=n-2kt+2\geq k$, we next locate an isolated vertex of the blue graph on $Z_v$. The maximality property above shows that, for each edge of $M_v$, at most one vertex of $X_v$ is incident with a blue edge from its endpoints that remain in $Z_v$. As $X_v$ itself is blue independent, at most $|M_v|\le k-1$ vertices of $X_v$ are nonisolated in $K_N[Z_v]$. Hence there exists $u \in X_v$ with $|N_B(u)\cap Z_v|=0$. Define $M_u$, $X_u$, $Y_u$, and $Z_u$ analogously. Since $Z_u\subseteq N_B(u)$, we have $Z_u\cap Z_v =\emptyset$. 
 Set $Z_1=Z_v$ and $Z_2=Z_u$. Then $Z_1$ and $Z_2$ are the asserted sets.
\end{proof}
Let $W=S\setminus(Z_1\cup Z_2)$. Clearly, $|W|\leq |S|-2(n-2kt+k+1)
  =2kt-3.$

 By Claims \ref{1} and \ref{Z1Z2}, 
for each $v \in Z_i$ and $i=1,2$, we have 
\begin{equation}
\begin{split}
  |N_B(v)\cap Z_{3-i}| &\geq |N_B(v)\cap S|-|N_B(v)\cap Z_i|-|W|\\
  &\geq (n-2kt+2k)-2(k-1)-(2kt-3)\\
  &\geq n-4kt+5.
  \label{eq:z degree_bound}
\end{split}
\end{equation}

If $W\not=\emptyset$, then for each $w \in W$, we have
\begin{equation}
\begin{split}
  |N_B(w)\cap(Z_1 \cup Z_2)| &\geq |N_B(w)\cap S|-(|W|-1)\\
  &\geq (n-2kt+2k)-(2kt-3-1)\\
  &=n-4kt+2k+4. 
  \label{eq:w degree_bound}
\end{split}
\end{equation}

Since $K_N[S]$ contains no blue $F_k$, no vertex $w \in W$ can simultaneously have many blue neighbors in both $Z_1$ and $Z_2$.
Thus we get the following claim.
 \begin{claim} \label{W0} For each $w \in W$ and $i=1,2$, if $|N_B(w)\cap Z_i|\geq 4kt+2k-7$, then $|N_B(w)\cap Z_{3-i}|\leq k-1$.   \end{claim}
 \begin{proof}         
 By symmetry of $Z_1$ and $Z_2$, we only need to prove the case for $i=1$. Suppose to the contrary that some $w \in W$ satisfies that $|N_B(w)\cap Z_2|\geq k$. Then for each $y \in N_B(w)\cap Z_2$, by Claim \ref{Z1Z2}
 and \eqref{eq:z degree_bound}, the blue neighbors of $y$ in $N_B(w)\cap Z_1$ is at least          \begin{align*} |N_B(y)\cap(N_B(w)\cap Z_1)| &\geq |N_B(y)\cap Z_1|+|N_B(w)\cap Z_1|-|Z_1|\\ &\geq (n-4kt+5)+(4kt+2k-7)-(n+k-2)\\&= k.\end{align*}          Thus we can greedily find a blue $kK_2$ between $N_B(w)\cap Z_1$ and $N_B(w)\cap Z_2$, which, together with $w$, forms a blue $F_k$, a contradiction.\end{proof}
By \eqref{eq:w degree_bound} and $n\geq 12kt+2k-19$, for each vertex $w\in W$, we have $|N_B(w)\cap Z_i|\geq 4kt+2k-7 \text{ for some $i\in \{1,2\}$}$. Combining Claim \ref{W0}, this implies that each $w\in W$ satisfies $|N_B(w)\cap Z_i|\leq k-1 \text{ for some $i\in \{1,2\}$}$. 
 Thus we can partition $W$ into two disjoint parts $W_1$ and $W_2$ by assigning each $w\in W$ to an index $i$ for which $|N_B(w)\cap Z_i|\leq k-1$ (choosing either index if both qualify). Consequently,
$$|N_B(w)\cap Z_i|\leq k-1\qquad\text{for every }w\in W_i.$$
Set $S_i=Z_i \cup W_i$ for $i=1,2$. We get the following claim, which implies that $S_1$ and $S_2$ induce a nearly blue complete bipartite graph.
\begin{claim} \label{S1S2}
    For each $v\in S_i$ and $i=1,2$, we have $|N_B(v)\cap S_{3-i}|\geq n-4kt+5$.
\end{claim}
\begin{proof}
For each $v\in S_i$ and $i=1,2$, if $v \in Z_i$, then
\begin{align*}
 |N_B(v)\cap S_{3-i}|&\geq |N_B(v)\cap S|-|N_B(v)\cap Z_i|-|W_i|\\
 &\geq n-2kt+2k-2(k-1)-(2kt-3)\\
 &=n-4kt+5.   
\end{align*}
If $v \in W_i$, then
\begin{align*}
    |N_B(v)\cap S_{3-i}|&\geq |N_B(v)\cap S|-(k-1)-(|W_i|-1)\\
    &\geq n-2kt+2k-(k-1)-(2kt-3-1)\\
    &=n-4kt+k+5.
    \end{align*}
    Thus the assertion follows.  
\end{proof}

We also have $|S_i|\le n-1$ for $i=1,2$. Indeed, at most $3(k-1)$ vertices of $Z_i$ are incident with a blue edge inside $Z_i$, and at most $|W_i|(k-1)$ vertices of $Z_i$ are incident with a blue edge from $W_i$. Since
$$|Z_i|-|W_i|(k-1)-3(k-1)\ge n-2k^2t+k+1>0,$$
some vertex of $Z_i$ has blue degree zero in $K_N[S_i]$. If $|S_i|\ge n$, this vertex is the center of a red $K_{1,n-1}$, a contradiction.

Recall that $|N_B(v)| \geq n+t-1$ for all $v\in V(K_N)$.
Thus
\begin{equation*}
    |N_B(a)\cap S| \geq n+t-1-(|A|-1)=n-2kt+2k+1 \text{ for any $a\in A$}.
\end{equation*}
Since $K_N$ has no $tF_k$, each $F_k$ in $A$ cannot form a blue $2F_k$ with some vertices in $S$. Thus we have the following claim.
\begin{claim} \label{A2}
   For each $F_k$ in $A$, there is at most one vertex whose number of blue neighbors is at least $2k$ in both $S_1$ and $S_2$.
\end{claim}   
\begin{proof} 
Suppose that two vertices $a_1,a_2$ of the same copy of $F_k$ in $A$ each have at least $2k$ blue neighbors in both $S_1$ and $S_2$. Put $T=4kt+2k-6$. Since $|N_B(a_i)\cap S|\ge n-2kt+2k+1$ and $n\ge10kt+2k-13$, each $a_i$ has at least $T$ blue neighbors in one of $S_1,S_2$.

We first record the following consequence of Claim~\ref{S1S2}. If $|N_B(a_i)\cap S_j|\ge T$ and $Q\subseteq N_B(a_i)\cap S_{3-j}$ has order $k$, then every $q\in Q$ has at least
   \begin{equation}
\begin{split}
  & |N_B(a_i)\cap S_j| + |N_B(q)\cap S_j| - |S_j| \\
  \geq & (4kt+2k-6) + (n-4kt+5) - (n-1) 
  = 2k.
\end{split}
\label{eq:degree_bound}
\end{equation}

For each $a_i$, choose an index $j_i$ for which $|N_B(a_i)\cap S_{j_i}|\ge T$. If $j_1\ne j_2$, choose $k$ blue neighbors $Q_i$ of $a_i$ in $S_{3-j_i}$. The two sets lie in different parts. By \eqref{eq:degree_bound}, match $Q_1$ by blue edges into $(N_B(a_1)\cap S_{j_1})\setminus Q_2$ and match $Q_2$ into $(N_B(a_2)\cap S_{j_2})\setminus Q_1$. Each vertex to be matched has at least $k$ available candidates, so Hall's condition holds for each of these two matchings.

If $j_1=j_2$, assume by symmetry that both equal $1$. Choose disjoint $k$-sets $Q_i\subseteq N_B(a_i)\cap S_2$; this is possible because both blue degrees into $S_2$ are at least $2k$. Use \eqref{eq:degree_bound} first to match $Q_1$ into $N_B(a_1)\cap S_1$, and then to match $Q_2$ into $N_B(a_2)\cap S_1$ while avoiding the first set of matched vertices. Again every vertex has at least $k$ candidates at the moment it is matched.

In either case the two matchings, together with the centers $a_1,a_2$, form two vertex-disjoint blue fans. Together with the remaining $(t-2)F_k$ in $A$, they form a blue $tF_k$, a contradiction.
\end{proof}

For $i=1,2$,
we define 
$$ A_i=\{a\in A : |N_B(a)\cap S_i|\leq 2k-1\}.$$
Since for each $a\in A$ we have
 $|N_B(a)\cap S| \geq n-2kt+2k+1>4k-2$, the sets $A_1$ and $A_2$ are disjoint. Every vertex outside $A_1\cup A_2$ has at least $2k$ blue neighbors in both parts, so Claim \ref{A2} gives $|A_1|+|A_2| \geq |A|-(t-1)$.
 Moreover, $$|S_1 \cup A_1|+ |S_2 \cup A_2| \geq |S|+|A|-(t-1) \geq N-(t-1)=2n-1.$$
Thus either $|S_1 \cup A_1| \geq n$ or $|S_2 \cup A_2| \geq n$. By symmetry, assume that $|S_1 \cup A_1| \geq n$.

   Recall that $S_1=Z_1 \cup W_1$ and every vertex of $W_1$ has at most $k-1$ blue neighbors in $Z_1$. 
   Thus the number of vertices of $Z_1$ incident with a blue edge from $W_1\cup A_1$ is at most $|W_1|(k-1)+|A_1|(2k-1)$, and 
   \begin{align*}
  |W_1|(k-1)+|A_1|(2k-1)
  &\leq (2kt-3)(k-1)+(t-1)(2k+1)(2k-1).
   \end{align*} 
Moreover, by 
the definition of $Z_1$, there 
are at most $3(k-1)$ vertices whose blue degree is not equal to $0$ in $K_N[Z_1]$.
Since $n \geq 6k^2t-4k^2-k-t+1$, the number of vertices with blue degree $0$ in $K_N[S_1 \cup A_1]$ is at least
\begin{align*}
   &|Z_1|-(2kt-3)(k-1)-(t-1)(2k+1)(2k-1)-3(k-1)\\
   \geq &n-6k^2t+4k^2+k+t
   \geq 1.
    \end{align*}
   Then by $|S_1 \cup A_1| \geq n$, we can find a red $K_{1,n-1}$ in $K_N[S_1\cup A_1]$, a contradiction.

 The proof of Theorem \ref{star tfan} is complete. \qed

\section{Proof of Theorem \ref{mr2}}\label{section5}

\begin{lem} (Burr, Erd{\H{o}}s, and Spencer \cite{BES1975})
    Let $G$ and $H$ be graphs. For any positive integer $t$, we have 
    $r(G, tH)\leq r(G,H)+(t-1)|H|$.
\end{lem}
\begin{proof}
Consider a red-blue coloring of a complete graph on $r(G,H)+(t-1)|H|$ vertices, and choose a maximum family of $s$ vertex-disjoint blue copies of $H$. If $s\geq t$, then there is a blue $tH$. Otherwise, deleting these copies leaves at least $r(G,H)$ vertices and no blue $H$; hence the remaining graph contains a red $G$.
\end{proof}

Then, from the exact value of 
$r(G,F_k)$ given by Theorem \ref{mr} and the upper bound provided by Lemma \ref{upper}, we obtain the following corollaries.
\begin{cor} \label{1upper tF_k}
    Let $k$ be a positive integer and $n\geq 36k^4$. If $G$ is a connected graph with $n$ vertices
and at most $ n(1+1/(204k^3+126k^2))$ edges,
then
$$r(G,tF_k)\leq 2n-1+(t-1)(2k+1).$$
\end{cor}

\begin{cor} \label{3upper tF_k}
    If $G$ is a simple graph on $n\geq1$ vertices with $m$ edges, then for any positive integers $k,t$,
     $$r(G,tF_k) \leq n+2m k-\frac{2m}{n}+(t-1)(2k+1).$$
\end{cor}

Note that $tF_k\subseteq \overline{K}_{kt}+tK_{1,k}$, the following lemma can help us first identify a blue $tK_{1,k}$ and then find a blue $tF_k$. 
\begin{lem} (Huang, Zhang, and Chen \cite{HZC}) \label{t star}
Let $k,t\geq1$, let $g(k,t)=\max\{21tk-3k+6,24k-12\}$, and let $n\geq 28t^2k^3$. If $G$ is a connected graph with $n$ vertices and at most $n(1+1/g(k,t))$ edges, then
$r(G,tK_{1,k})\leq n+k+t-2$.  
\end{lem}

We now begin to  prove Theorem \ref{mr2}.

 Since $\chi(tF_k)=3$ and $s(tF_k)=t$, \eqref{lower} gives the matching lower bound. It therefore suffices to prove the upper bound. The proof will be by induction on $t$. When $t=1$, it is verified
 by Theorem \ref{mr}. Assuming the theorem holds for $t -1$, we now proceed to consider the case where $k\geq 1$ and $t \geq 2$.

 To apply the inductive method, it is necessary to ensure that the lower bound on $n$ and the upper bound on $e(G)$ are
 compatible with induction. This can be derived from the facts that the lower bound on $n$ is increasing in $t$ and the upper bound on $e(G)$ is decreasing in $t$.
 
 Let $N=2n+t-2$. Given a red-blue edge-coloring of $K_N$, assume that $K_N$ contains neither a red $G$ nor a blue $tF_k$. 
 By Theorem \ref{star tfan}, the case where $G$ is a star is already settled. Thus we assume $G$ is not a star, and then 
 by Lemma \ref{trichotomy}, our argument divides into three cases.

\vskip 2mm
\noindent {\textbf{Case 1.}} $G$ has a suspended path with at least $(2k^2+3k+1)t-1+(kt-k+1)$ vertices.
\vskip 2mm
 Denote this suspended path by $P$.
 Let $H$ be the graph obtained from $G$ by shortening $P$ by $kt-k+1$ vertices. Since $n\geq 161t^2k^4$ and $t\geq 2$, we have that
 \begin{align*}
 |H| &=n-(kt-k+1) \geq 36k^4 ,\\
   e(H) 
  &\leq n\left(1+\frac{1}{204 t k^3 + 147 t k^2 }\right)-(kt-k+1)  \\
   &\leq \left(n-(kt-k+1) \right) \left(1+\frac{1}{204k^3+126k^2}\right).
 \end{align*}
 Then by Corollary \ref{1upper tF_k}, $$r(H,tF_k)\leq 2(n-(kt-k+1))+(t-1)(2k+1)= 2n+t-3<N.$$ Therefore, $K_N$ contains a red $H$. Let $H'$ be the subgraph of $K_N[R]$ which can
 be obtained from $H$ by lengthening the suspended path $P$ as much as possible (up to
 $kt-k+1$ vertices). And let $P'$ be the suspended path in $H'$. If $H'$ is $G$, the proof of this case is complete. If not, there are two possibilities to consider, according to whether $K_N-H'$ has a blue $tK_{1,k}$ or not.
 
 If $K_N-H'$ contains a blue $tK_{1,k}$, write these vertex-disjoint stars as $Y_1,\ldots,Y_t$, each of order $k+1$. The red path $P'$ has at least
 $t(2k^2+3k)+(t-1)$
 vertices. We may therefore choose, in their natural order along $P'$, pairwise vertex-disjoint subpaths $Q_1,\ldots,Q_t$, each of order $(2k+3)k=2k^2+3k$, with one unused path vertex between consecutive subpaths. Apply Lemma~\ref{path extension} to $Q_i\cup Y_i$ with $b=k+1$, $c=2k+1$, and $d=k$. Its first alternative replaces $Q_i$ by a red path with the same endvertices and one additional vertex, contradicting the maximality of $P'$. Its second alternative gives a blue $K_{2k+1}$, and its third gives a blue $\overline K_k+K_{1,k}$; either graph contains a blue $F_k$. These $t$ fans lie in the pairwise disjoint sets $Q_i\cup Y_i$, yielding a blue $tF_k$, a contradiction.

 Now we consider the case that $K_N-H'$ contains no blue $tK_{1,k}$. Let $S=V(K_N-H')$. Because $ |H'| \le n-1$, we have $|S|\geq N-(n-1)= n+t-1$.

Let $H_0$ be the graph obtained from $G$ by shortening $P$ by $k-1$ vertices. Since $n\ge 161t^2k^4$, we have \begin{align*}  
|H_0| &=n-k+1 \geq 28t^2k^3 ,\\       
e(H_0) 
&\leq n\left(1+\frac{1}{204 t k^3 + 147 t k^2}\right)-(k-1) 
 \leq \left(n-k+1 \right) \left(1+\frac{1}{21tk-3k+6}\right).  \end{align*}         
Then by Lemma \ref{t star}, $$r(H_0,tK_{1,k})\leq (n-k+1)+k+t-2=n+t-1\le |S|.$$ Because $K_N[S]$ contains no blue $tK_{1,k}$, there must be a red $H_0$ in $K_N[S]$. 
If $k=1$, then $H_0=G$, which is already a contradiction. Hence assume $k\geq2$ in the remainder of this subcase.
Let $H_0'$ be the subgraph of $K_N[S]$ obtained from $H_0$ by lengthening the shortened suspended path through red edges as many as possible (up to $k-1$ vertices). Since $K_N$ has no red $G$, $H_0'\not=G$. Let $P_0'$ be the suspended path of $H_0'$. Since $|S|-|H_0'| \geq n+t-1-(n-1)=t$, choose distinct vertices $v_1,\ldots,v_t\in S\backslash V(H_0')$. By the maximality of $P_0'$, no $v_i$ is red-adjacent to two consecutive vertices of $P_0'$, and hence each $v_i$ has at least $\lfloor|V(P_0')|/2\rfloor\ge kt$ blue neighbors on $P_0'$. Choose $k$ such neighbors for each center successively; before the last choice, at most $(t-1)k$ path vertices have been used, so at least $k$ candidates remain. This gives a blue $tK_{1,k}$ in $K_N[S]$, a contradiction.

 \vskip 2mm
  \noindent  {\textbf{Case 2.}} $G$ has a matching consisting of $3kt-k-2$ end-edges.
 \vskip 2mm
By the induction hypothesis, we have
 $r(G,(t-1)F_k)\leq 2n+t-3<N$,
 which implies that the graph $K_N$ contains a blue subgraph $(t-1)F_k$. We denote the
 vertex set of this subgraph by $A$ and $V(K_N)-A$ by $S$.

Let $H$ be the graph obtained from $G$ by deleting $3kt-k-2$ end-vertices of this matching. Since $n\geq 161t^2k^4$ and $t\geq 2$, we have that
 \begin{align*}
 |H| &=n-(3kt-k-2) \geq 36k^4,\\
   e(H) 
   &\leq n \left(1+\frac{1}{204 t k^3 + 147 t k^2 }\right)-(3kt-k-2)  \\
   &\leq (n-(3kt-k-2))\left(1+\frac{1}{204k^3+126k^2}\right).
 \end{align*}
By Theorem \ref{mr}, we get 
\begin{align*}
    r(H,F_k)&\leq 2(n-(3kt-k-2))-1 \leq 2n-2kt+2k-1= |S|.
    \end{align*}
Thus $K_N[S]$ contains a red $H$. 
Let $X$ be the set of $3kt-k-2$ vertices of $H$ which were adjacent to the end-vertices deleted from $G$ and $Y=V(K_N)-V(H)$. By Lemma \ref{matching}, there is either a red matching that covers $X$, or a blue $K_{c+1,|Y|-c}$, where $0 \leq c \leq |X|-1$. 
If the former holds, then $K_N$ would contain a red $G$, leading to a contradiction. If the latter holds, let $X'$ be the vertices of $X$ and $Y'$ be the vertices of $Y$ contained in the blue $K_{c+1,|Y|-c}$.
 Since $K_N[S]$ contains no blue $F_k$, there is no blue $kK_2$ in  $K_N[N_B(v)\cap S]$ for any $v\in S$, and then
 $|N_B(v)\cap S|\leq n+k-2$ by Lemmas \ref{K2} and \ref{kK2}. Thus $|Y'|=|Y|-c\leq n+k-2+|A|$.
 Then \begin{align*}
     c&\geq |Y|-n-k+2-|A|\\
     &=2n+t-2-(n-|X|)-n-k+2-(t-1)(2k+1)\\
     &= |X|+t-k-(t-1)(2k+1)\\
     &= 3kt-k-2+t-k-(t-1)(2k+1)\\
     &= kt-1.
     \end{align*} Moreover, 
 $$|Y'|=|Y|-c\geq 2n+t-2-(n-|X|)-(|X|-1)=n+t-1.$$
 If $K_N[Y']$ contains a blue $ktK_2$, then by $|X'|=c+1\geq kt\geq t$, there is a blue $\overline{K}_{t}+ktK_2$. And if $K_N[Y']$ contains a blue $tK_{1,k}$, then by $|X'|=c+1\geq kt$, there is a blue $\overline{K}_{kt}+tK_{1,k}$. In both cases, $K_N$ contains a blue $tF_k$. Thus $K_N[Y']$ contains neither a blue $ktK_2$ nor a blue $tK_{1,k}$. 
 
 Next we claim that $K_N[ Y']$ contains a red $G$.
Let $M$ be the maximum blue matching in $K_N[Y']$ and $Y''=Y'-V(M)$. Since $K_N[Y']$ contains no blue $ktK_2$, $|M|=s\leq kt-1$ and $K_N[Y'']$ is red complete. 
Let $M = \{x_1y_1,x_2y_2,\ldots,x_sy_s\}$ and assume $|N_B(x_i)\cap Y''|\le |N_B(y_i)\cap Y''|$ for all $1\le i\le s$. Then by the maximality of $M$, either $|N_B(x_i)\cap Y''|= |N_B(y_i)\cap Y''|=1$ and $x_i$ and $y_i$ are
 adjacent to the same vertex in $Y''$, or $|N_B(x_i)\cap Y''|=0$. Set
 $$U=\{y_i: |N_B(y_i)\cap Y''| \geq kt\},\qquad u=|U|.$$
If $u\ge t$, any $t$ vertices of $U$ are centers of vertex-disjoint blue $K_{1,k}$'s with leaves in $Y''$: for every set of $r\le t$ centers, the union of their blue neighborhoods has order at least $kt\ge rk$, so the required replicated form of Hall's condition holds. Therefore $u\le t-1$. Choose any $U'\subseteq Y''$ with $|U'|=u$. We claim that there is a red matching from $V(M)\setminus U$ into $Y''\setminus U'$ covering $V(M)\setminus U$.

In fact, for each $v\in V(M)\setminus U$, we have
\begin{align*}
    |N_R(v)\cap (Y''\setminus U')| &\geq |Y''|-u-(kt-1)\\
    &=|Y'|-2s-u-kt+1\\
    &\geq n+t-1-2(kt-1)-(t-1)-kt+1\\
    &=n-3kt+3\\
    &\geq 2kt-2\\
    &\geq |V(M)\setminus U|.
\end{align*}
 Thus we can greedily find such a red matching; choose its target vertices to be distinct. Since $K_N[Y'']$ is red complete, we may also match every vertex of $U'$ by a red edge to a new, distinct vertex of $Y''\setminus U'$, avoiding the targets already chosen. There are enough such new targets because
 $$|Y''|-2s=|Y'|-4s\geq n+t-1-4(kt-1)\geq t-1\geq u.$$

Let $H'$ be obtained from $G$ by deleting $2s$ end-vertices from the fixed matching of $3kt-k-2$ end-edges. This is possible because $3kt-k-2\geq 2kt-2\geq 2s$. The deleted leaves have $2s$ distinct support vertices. Also,
$$|Y''\setminus U'|=|Y'|-2s-u\geq n+t-1-2s-u\geq n-2s=|H'|.$$
Embed $H'$ into the red clique $Y''\setminus U'$ so that its support vertices are mapped to the $2s$ chosen matching targets. Map $2s-u$ deleted leaves to $V(M)\setminus U$ and the remaining $u$ leaves to $U'$. All required incident edges are red, so this is a red copy of $G$ in $K_N[Y']$, a contradiction.

 \vskip 2mm
\noindent {\textbf{Case 3.}} 
$G$ has neither a suspended path of $(2k^2+3k+1)t-1+(kt-k+1)$ vertices, nor a matching formed by $3kt-k-2$ end-edges.
\vskip 2mm 
For convenience, we let $e(G)=n+\ell$, $q=(2k^2+3k+1)t-1+(kt-k+1)=2k^2t+4kt+t-k$ and $s=3kt-k-2$, where $\ell \leq n/(204 t k^3 + 147 t k^2)$.

After removing all vertices of degree 1 and their incident edges from the graph $G$, we
 obtain a graph denoted by $G'$. Clearly, $G'$ is a connected graph. Because $n\geq 161t^2k^4$,
 then by Lemma \ref{trichotomy},
\begin{align*}
   |G'| &\leq \gamma=(q-2)(2s+3\ell-2)+1\\
   &\leq (q-2)(2s-2)+1+\frac{3(q-2)}{204 t k^3 + 147 t k^2}n \\
   &\leq \frac{96k^2t}{204 t k^3 + 147 t k^2}n+\frac{3(2k^2t+4kt+t-k-2)}{204 t k^3 + 147 t k^2}n\\
   & = \frac{102k^2t+12kt+3t-3k-6}{204 t k^3 + 147 t k^2}n.
\end{align*}
Moreover,
 \begin{align*}
 |G'|-1 \leq e(G') &\leq |G'|+\frac{n}{204 t k^3 + 147 t k^2}
 \leq \frac{102k^2t+12kt+3t-3k-5}{204 t k^3 + 147 t k^2}n . 
\end{align*}
By Lemma \ref{trichotomy},  $G$ has a vertex adjacent to at least 
$$\left\lceil \frac{n-\gamma}{3kt-k-3} \right\rceil$$ 
vertices of degree 1.
Denote this vertex by $v$. By Theorem \ref{star tfan}, there exists a red star $K_{1,n-1}$ in $K_N$, with the center
 vertex denoted by $x$. By Corollary \ref{3upper tF_k}, if $|G'| \geq 2$, then
\begin{align*}
    r(G',tF_k)&\leq |G'|+2ke(G')-\frac{2e(G')}{|G'|}+(t-1)(2k+1)\\
    &\leq |G'|+2ke(G')-1+\frac{24}{204 t k^3 + 147 t k^2}n\\
    &\leq \frac{204k^3t + 126k^2t - 6k^2 + 18kt - 13k + 3t + 18}{204 t k^3 + 147 t k^2}n-1\\
    &\leq  n-1.
\end{align*}
Thus $K_N[N_R(x)]$ contains a red $G'$. Moreover, this holds trivially if $|G'| = 1$. Replacing $v$ in this $G'$ by $x$ gives another copy of red $G'$. Every vertex of $G''-G'$ is an original leaf of $G$ whose support lies in $G'$. Embed these leaves one at a time, always after their support has been embedded. In this way we greedily embed $G''$ into
 the red subgraph of $K_N$, where $G''$ is obtained by removing all degree 1 vertices adjacent to $v$ in the graph $G$. Clearly, $G''$ contains $G'$ as a subgraph and
 $$|G''| \leq n-\left\lceil \frac{n-\gamma}{3kt-k-3}\right\rceil.$$
If at some step the embedding cannot proceed, the image $z$ of the relevant support vertex has no unused red neighbor and hence is adjacent to at most $|G''|-2$ vertices by red edges. In other words, the
 number of blue edges incident to vertex $z$ is at least $N-|G''|+1$. Because $n\geq 161t^2k^4$, we have
 $n-\gamma \ge (3kt-k-3)(2kt-k-1)$. Then
 \begin{align*}
    N-|G''|+1 &\geq 2n+t-2-n+\left\lceil \frac{n-\gamma}{3kt-k-3}\right\rceil +1\\
    &=n+t-1+\left\lceil \frac{n-\gamma}{3kt-k-3}\right\rceil \\
    & \geq n+t-1+(2kt-k-1)\\
    & = n+k-1+(t-1)(2k+1)\\
    & = r(G,kK_2)+|(t-1)F_k|.
 \end{align*}
By the induction hypothesis, we have
 $r(G,(t-1)F_k)\leq 2n+t-3=N-1$,
 which implies that the graph $K_N\setminus \{z\}$ contains a blue subgraph $(t-1)F_k$. We denote the
 vertex set of this subgraph by $A$. In the blue neighborhood of vertex
 $z$, there are at least $n+k-1$ vertices that do not belong to the set $A$. By Lemmas \ref{K2} and \ref{kK2}, in the graph induced by these vertices, there exists either a red subgraph
 isomorphic to $G$, or a blue $kK_2$, which, together with vertex $z$ and the vertices in $A$,
 induces a blue subgraph containing $tF_k$. Both cases lead to a contradiction.

 Thus, $G''$ can always be embedded into the red subgraph of $K_N$. Next, we only need
 to embed the degree-1 vertices adjacent to $v$, and the entire graph $G$ can be embedded
 into the red subgraph of $K_N$. This is feasible because the vertex $x$ into which $v$ is
 embedded is adjacent to at least $n-1$ red edges, so there are enough vertices to embed the degree-1 vertices of $v$. This completes the proof of Theorem \ref{mr2}. \qed

\section{Concluding Remark}
The numerical hypotheses in Theorems~\ref{mr} and~\ref{mr2} arise from the third case of each proof, where the Trichotomy Lemma is combined with the general Ramsey upper bound of Lemma~\ref{upper}. The constants have been chosen so that every displayed estimate is uniform over the full stated parameter range. They are not expected to be optimal; substantial improvements would require either a sharper structural lemma for sparse graphs or a stronger general upper bound for Ramsey numbers against fans.

 \section*{Acknowledgments}
This research was supported  by National Key R\&D Program of China under grant number 2024YFA1013900 and NSFC under grant number 12471327.


\begin{thebibliography}{99}
\small \setlength{\itemsep}{-.10mm}

\bibitem{AKS1980}
M. Ajtai, J. Koml{\'o}s, and E. Szemer{\'e}di, A note on Ramsey numbers, J. Combin. Theory Ser. A 29 (1980), 354--360 . 
\bibitem{AKS1981}
M. Ajtai, J. Koml{\'o}s, and E. Szemer{\'e}di, A dense infinite Sidon sequence, European J. Combin. 2 (1981), 1--11.



\bibitem{BE}
J.A. Bondy and P. Erd\H{o}s, Ramsey numbers for cycles in graphs, J. Combin. Theory
 Ser. B 14 (1973), 46--54.
\bibitem{Brennan}
M. Brennan, Ramsey numbers of trees and unicyclic graphs versus fans, Discrete Math. 340 (2017), 969--983.
 
\bibitem{Burr1981}
S.A. Burr, Ramsey numbers involving graphs with long suspended paths,
J. London Math. Soc. (2) 24 (1981), 405--413.


\bibitem{BEFRS}
S.A. Burr, P. Erd\H{o}s, R.J. Faudree, C.C. Rousseau, and R.H. Schelp, Ramsey numbers for the
pair sparse graph-path or cycle, Trans. Amer. Math. Soc. 269 (1982), 501--512.

\bibitem{BES1975}
S.A. Burr, P. Erd\H{o}s, and  J. H. Spencer, Ramsey theorems for multiple copies of graphs, Trans. Amer. Math. Soc. 209 (1975), 87--99.

\bibitem{CJMS2025}
M. Campos, M. Jenssen, M. Michelen, and J. Sahasrabudhe, A new lower bound for the Ramsey numbers $R(3,k)$, arXiv:2505.13371 (2025).




\bibitem{CL2025}
F. Chung and Q. Lin, Fan-complete Ramsey numbers, Adv. in Appl. Math. 171 (2025), 102939.

\bibitem{Chvatal1977}
V. Chv\'atal, 
Tree-complete graph Ramsey numbers, 
J. Graph Theory 1 (1977), 93.

\bibitem{Chvatal1972}
 V. Chv\'atal and F. Harary, Generalized Ramsey theory for graphs. III. Small off
diagonal numbers, Pacific J. Math. 41 (1972), 335--345.


\bibitem{EFRS1985}
P. Erd\H{o}s, R.J. Faudree, C.C. Rousseau, and R.H. Schelp,
Multipartite graph-sparse graph Ramsey numbers,
Combinatorica 5 (1985), 311--318.

\bibitem{FSS}
R.J. Faudree, R.H. Schelp, and J. Sheehan, Ramsey numbers for matchings, Discrete Math. 32 (1980), 105--123.



\bibitem{Hall}
P. Hall, 
On representatives of subsets, 
 J. London Math. Soc. (1) 10 (1935), 26--30.
\bibitem{HP2023}
S. Hu and Y. Peng,
 Ramsey numbers of stripes versus trees and unicyclic graphs,
J. Oper. Res. Soc. China 13 (2025), 297--312. 

\bibitem{HZC}
T. Huang, Y. Zhang, and Y. Chen, Minimum degree and sparse connected spanning subgraphs, arXiv:2507.03264v2 (2025).

\bibitem{KOS2019}
S. Kadota, T. Onozuka, and Y. Suzuki, The graph Ramsey number $R(F_l,K_6)$, Discrete Math. 342 (2019), 
1028--1037.

\bibitem{Kim1995}
J. H. Kim, The Ramsey number $R(3, t)$ has order of magnitude $t^2/ \log t$, Random Struct. Alg. 7 (1995), 173--207.


\bibitem{LR1996}
Y. Li and C.C. Rousseau, Fan-complete graph Ramsey numbers. J. Graph Theory 23 (1996), 413--420.




\bibitem{S1983}
J. Shearer, A note on the independence number of triangle-free graphs, Discrete 
Math. 46 (1983), 83--87.


\bibitem{Shi 2010}
L. Shi, Ramsey numbers of long cycles versus books or wheels, European J. Combin. 31 (2010), 828--838.

\bibitem{BB2005}
 Surahmat, E.T. Baskoro, and H.J. Broersma, The Ramsey numbers of fans versus $K_4$, Bull. Inst. Combin. Appl. 43 (2005), 96--102.

\bibitem{WQ2022}
M. Wang and J. Qian, Ramsey numbers for complete graphs versus generalized fans, Graphs Comb. 38(6) (2022), 186.



\bibitem{ZBC}
Y. Zhang, H. Broersma, and Y. Chen, Ramsey numbers of trees versus fans, Discrete Math. 338 (6) (2015), 994--999.

\bibitem{ZC2023}
Y. Zhang and Y. Chen, Ramsey goodness of fans, Discrete Math. 349 (2026), 114771.





\bibitem{ZC}
Y. Zhang and Y. Chen, Trichotomy and $tK_m$-goodness of sparse graphs, J. Graph Theory 0 (2026), 1--13.
\url{https://doi.org/10.1002/jgt.70125}

\end{thebibliography}
\end{document}